\documentclass[12pt]{amsart}
\usepackage{amsfonts}
\usepackage{amssymb}
\usepackage{amsmath}
\input xy
\xyoption {all}

\usepackage{geometry}
\usepackage{amsthm}
\usepackage{amssymb}
\usepackage{latexsym}
\usepackage[dvips]{graphics}

\newtheorem{thm}{Theorem}[section] 
\newtheorem{cor}[thm]{Corollary}

\newtheorem{prop}[thm]{Proposition}
\newtheorem{remark}[thm]{Remark}

\newtheorem*{ack}{Acknowledgements}

\numberwithin{equation}{section}

\newcommand{\HC}{{\mathcal H}}
\newcommand{\VV}{{\mathcal V}}
\newcommand{\FC}{{\mathcal F}}

\newcommand{\Z}{\mathbb{Z}}

\newcommand{\Q}{\mathbb{Q}}
\newcommand{\R}{\mathbb{R}}
\newcommand{\C}{\mathbb{C}}

\newcommand{\HB}{\mathbb{H}}

\begin{document}
\date{\today}
\title[Euler characteristics]{Euler characteristics of algebraic varieties }

\author[Sylvain E. Cappell ]{Sylvain E. Cappell}
\address{S. E. Cappell: Courant Institute, New York University, New York, NY-10012}
\email {cappell@cims.nyu.edu}

\author[Laurentiu Maxim ]{Laurentiu Maxim}
\address{L. Maxim: Courant Institute, New York University, New York, NY-10012}
\email {maxim@cims.nyu.edu}

\author[Julius L. Shaneson ]{Julius L. Shaneson}
\address{J. L. Shaneson: Department of Mathematics, University of Pennsylvania, Philadelphia, PA-19104}
\email {shaneson@sas.upenn.edu}

\thanks{The first and third authors partially supported by grants from NSF and DARPA.
The second author partially supported by a grant from the NYU
Research Challenge Fund.}

\begin{abstract}
This note studies the behavior of Euler characteristics and of
intersection homology Euler characterstics under proper morphisms of
algebraic (or analytic) varieties. The methods also yield, for
algebraic (or analytic) varieties, formulae comparing these two
kinds of Euler characteristics. The main results are direct
consequences of the calculus of constructible functions and
Grothendieck groups of constructible sheaves. Similar formulae for
Hodge theoretic invariants of algebraic varieties under morphisms
were announced by the first and third authors in \cite{CS1, S}.
\end{abstract}

\maketitle

\section{Introduction}

We study the behavior of (intersection homology) Euler
characteristics under proper morphisms of complex algebraic (or
analytic) varieties. We begin by discussing simple formulae for the
usual Euler-Poincar\'{e} characteristic, then show that similar
formulae hold for the intersection homology Euler characteristic, as
well as for the corresponding Chern homology classes of MacPherson.
The methods used in the present paper also yield formulae expressing
the Euler characteristics of usual and intersection homology of an
algebraic (or analytic) variety in terms of each other and
corresponding invariants of the subvarieties formed by the closures
of its singular strata.

The main results of this note are direct applications of the
standard calculus of constructible functions and Grothendieck groups
of constructible sheaves. Some of the formulae on the intersection
homology Euler characteristic were originally proven (cf.
\cite{CMS0})  with the aid of the deep BBDG decomposition theorem
for the pushforward of an intersection homology complex under a
proper morphism (cf. \cite{BBD,Cat}). The functorial approach
employed here was suggested by the referee. However, the core
calculations used in proving these results are modeled on our
original approach based on BBDG.

This note is a first step in an ongoing project that deals with the
study of genera of complex algebraic (or analytic) varieties. In the
forthcoming papers \cite{CMS} and \cite{CLMS} we will discuss the
behavior of Hodge theoretic genera under proper morphisms, and
provide explicit formulae for the pushforward of various
characteristic classes. The functorial approach and the language of
Grothendieck groups of constructible sheaves used in this paper
allow a simple translation of the underlying ideas to the
forthcoming papers, where Grothendieck groups of Saito's algebraic
mixed Hodge modules will be employed.

Unless otherwise specified, all homology and intersection homology
groups in this paper are those with rational coefficients. We assume
the reader's  familiarity with intersection homology, and for some
arguments also with (Grothendieck groups of) constructible sheaves
and derived categories. However, our results are also explained in
the simpler language of constructible functions, which only relies
on Euler characteristic information.

\begin{ack}We are grateful to the anonymous referee for his valuable
comments and suggestions regarding the functorial approach used in
the present paper. \end{ack}

\section{Topological Euler-Poincar\'{e} characteristic}
For a complex algebraic variety $X$, let $\chi(X)$ denote its
topological Euler characteristic. Then $\chi(X)$  equals  the
compactly supported Euler characteristic, $\chi_c(X)$ (cf.
[\cite{F}, p.141], [\cite{Sc}, \S6.0.6]).\footnote{This fact is not
true outside the category of complex varieties (e.g., if $X$ is an
oriented $n$-dimensional topological manifold, then Poincar\'{e}
duality yields that $\chi_c(X)=(-1)^n \chi(X)$).} The additivity
property for the Euler characteristic reads as follows: for  $Z$ a
Zariski closed subset of $X$, the long exact sequence of the
compactly supported cohomology
$$\cdots \to H^i_c(X \setminus Z) \to H^i_c(X) \to H^i_c(Z) \to H^{i+1}_c(X \setminus Z) \to \cdots$$
yields that $\chi_c(X)=\chi_c(Z)+\chi_c(X \setminus Z)$, therefore
the same relation  holds for $\chi$. The multiplicativity property
for fibrations asserts that if $F \hookrightarrow E \to B$ is a
locally trivial topological fibration such that  the three Euler
characteristics $\chi(B)$, $\chi(F)$ and $\chi(E)$ are defined, then
$\chi(E)=\chi(B)\cdot \chi(F)$ (e.g., see \cite{Di}, Corollary
2.5.5). In particular, if $f:X \to Y$ is a proper smooth submersion
of smooth manifolds, with $Y$ connected and generic fiber $F$, then:
\begin{equation} \chi(X)=\chi(Y) \cdot \chi(F).\end{equation}
Indeed, by Ehresmann's theorem, such a map is a locally trivial
fibration in the complex topology.

In this note we generalize this multiplicative property of proper
smooth submersions in two different directions: first, we study the
behavior of the usual Euler characteristic under arbitrary proper
maps of possibly singular varieties; second, we replace the usual
cohomology by intersection cohomology when dealing with singular
varieties, and study the behavior of the intersection homology Euler
characteristic under arbitrary proper morphisms. The formulae we
obtain here are classically referred to as the \emph{stratified
multiplicative property} for Euler characteristics (cf.
\cite{CS1,S}).

\bigskip

Let $Y$ be a topological space with a finite partition $\VV$ into a
disjoint union of finitely many connected subsets $V$ satisfying the
\emph{frontier condition}: ``$V \cap \bar{W} \neq \emptyset$ implies
that $V \subset \bar{W}$". (The main examples of such spaces are
complex algebraic or compact analytic varieties with a fixed Whitney
stratification.) Then $\VV$ is partially ordered by ``$V \leq W$ if
and only if $V \subset \bar{W}$". Let $F_{\VV}(Y)$ be the abelian
group of $\VV$-constructible functions on $Y$, i.e., of functions
$\alpha:Y \to \Z$ such that $\alpha|_V$ is constant for all $V \in
\VV$. This is a free abelian group with basis $\{1_V | V \in \VV
\}$, so that
$$\alpha=\sum_{V \in \VV} \alpha(V) \cdot 1_V.$$
Note that $\{ 1_{\bar V} | V \in \VV \}$ is another basis for
$F_{\VV}(Y)$, since $$1_{\bar V}= \sum_{W \leq V} 1_W$$ and the
matrix $A=(a_{W,V})$, with $a_{W,V}:=1$ for $W \leq V$ and $0$
otherwise, is upper triangular with respect to $\leq$, with all
diagonal entries equal to $1$. Thus $A$ is invertible. The non-zero
entries of $A^{-1}=(a'_{W,V})$ can inductively be calculated (e.g.,
see \cite{Sta}, Prop. 3.6.2) by $a'_{V,V}=1$ and, for $W <V$,
$$a'_{W,V}=-\sum_{W \leq S < V} a'_{W,S} \cdot a_{S,V}.$$
This implies the following
\begin{prop}\label{one} For each $V \in \VV$, define inductively
$\hat{1}_{\bar V}$ by the formula $$\hat{1}_{\bar V}=1_{\bar V} -
\sum_{W<V} \hat{1}_{\bar W}.$$ Then, for any $\alpha \in
F_{\VV}(Y)$, one has the equality
\begin{equation}\label{eq1} \alpha=\sum_V \alpha(V) \cdot \hat{1}_{\bar
V}.\end{equation}
\end{prop}
\begin{proof}
As the notation indicates, $\hat{1}_{\bar V}$ depends only on the
space $\bar V$ with its induced partition. Then by the above
considerations we have $$\alpha=\sum_V \alpha(V) \cdot 1_V=\sum _{W
\leq V} \alpha(V) \cdot a'_{W,V} \cdot 1_{\bar W},$$ and formula
(\ref{eq1}) follows from the inductive identification (for V fixed):
$$\sum_{W \leq V} a'_{W,V} \cdot 1_{\bar W}=1_{\bar V}-\sum_{W \leq
S <V} a'_{W,S} \cdot a_{S,V} \cdot 1_{\bar W} = \hat{1}_{\bar V}.$$

\end{proof}

\begin{remark}\rm
\medskip \noindent$(1)$ If there is a stratum $S \in \VV$ which is dense in $Y$, i.e.,
$\bar S=Y$, or $V \leq S$ for all $V \in \VV$, then formula
(\ref{eq1}) can be rewritten as
\begin{equation}\label{eq2}
\alpha = \alpha(S) \cdot 1_Y + \sum_{V <S} \left( \alpha(V) -
\alpha(S) \right) \cdot \hat{1}_{\bar V}.
\end{equation}

\medskip \noindent$(2)$ For a group homomorphism $\phi :
F_{\VV}(Y) \to G$ for some abelian group $G$, one obtains similar
descriptions for $\phi(\alpha)$ in terms of
$$\hat{\phi}(\bar V):=\phi(\hat{1}_{\bar V})=\phi(1_{\bar V}) -
\sum_{W<V} \phi(\hat{1}_{\bar W}).$$

\end{remark}

\bigskip

For the rest of this section we specialize to the complex algebraic
(or compact analytic context), with $Y$ a reduced complex algebraic
variety (or a reduced compact complex analytic space), and all $V
\in \VV$ locally closed constructible subsets. The group $F_c(Y)$ of
all complex algebraically (resp. analytically) constructible
functions is defined as the direct limit of these $F_{\VV}(Y)$. Then
one has the following important group homomorphisms on $F_c(Y)$
(e.g., see \cite{KS,M,Sc,SY}):
\begin{enumerate}
\item The Euler characteristic with compact support $\chi_c:F_c(Y)
\to \Z$ characterized by $\chi_c(1_Z)=\chi_c(Z)$ for $Z \subset Y$ a
locally closed constructible subset.
\item The Euler characteristic $\chi:F_c(Y) \to \Z$ characterized by
$\chi(1_Z)=\chi(Z)$ for $Z \subset Y$ a closed algebraic (resp.
analytic) subset.
\item For $f:X \to Y$ a proper complex algebraic (resp. analytic)
map, the functorial pushdown $f_*:F_c(X) \to F_c(Y)$ is
characterized by $f_*(1_Z)(y)=\chi(Z \cap \{f=y\})$ for $Z \subset
X$ a closed algebraic (resp. analytic) subset.
\item The Chern-MacPherson class transformation $c_*:F_c(Y) \to
H^{BM}_{2*}(Y;\Z)$, which commutes with proper pushdowns, and is
uniquely characterized by this property together with the
normalization $c_*(1_M)=c^*(TM) \cap [M]$, for $M$ a complex
algebraic (resp. analytic) manifold.
\end{enumerate}
In fact, as already pointed out, in this context we have
$\chi=\chi_c$. Moreover, $\chi \circ f_*=\chi$, and for $Y$ compact
one gets by functoriality  $\chi(\alpha)=\deg(c_*(\alpha))$, for any
$\alpha \in F_c(Y)$.

\bigskip

Now let $f:X \to Y$ be a proper complex algebraic (resp. analytic)
map, with $Y$ as above. Assume $f_*(\alpha) \in F_{\VV}(Y)$ for a
given $\alpha \in F_c(X)$ (e.g., $\VV'$ and $\VV$ are complex
Whitney stratifications of $X$ and resp. $Y$, such that $f$ is a
stratified submersion, and $\alpha \in F_{\VV'}(X)$). Then
$$f_*(\alpha)=\sum_{V \in \VV} f_*(\alpha)(V) \cdot 1_V,$$ with
$$f_*(\alpha)(V)=\chi(\alpha|_{F_V}),$$ for $F_V$ the fiber of $f$
over a point in $V$. Assume, moreover, that $Y$ is irreducible, so
there is a dense stratum $S \in \VV$, with $F:=F_S$ a general fiber
of $f$. In terms of $f_*(\alpha)$, formula (\ref{eq2}) yields the
following:
\begin{equation}\label{eq3}
f_*(\alpha)=\chi(\alpha|_F) \cdot 1_Y + \sum_{V < S} \left(
\chi(\alpha|_{F_V})-\chi(\alpha|_F) \right) \cdot \hat{1}_{\bar V}.
\end{equation}
By applying the homomorphism $\chi$ and resp. $c_*$ to the equation
(\ref{eq3}), we obtain the following formulae:
\begin{cor}
\begin{equation}\label{eq4}
\chi(\alpha)=\chi(\alpha|_F) \cdot \chi(Y) + \sum_{V < S} \left(
\chi(\alpha|_{F_V})-\chi(\alpha|_F) \right) \cdot \hat{\chi}({\bar
V}).
\end{equation}
\begin{equation}\label{eq5}
f_*(c_*(\alpha))=\chi(\alpha|_F) \cdot c_*(Y) + \sum_{V < S} \left(
\chi(\alpha|_{F_V})-\chi(\alpha|_F) \right) \cdot \hat{c}_*({\bar
V}),
\end{equation}
where $c_*(Y):=c_*(1_Y)$, and similarly for $\hat{c}_*({\bar V})$,
which by the functoriality of $c_*$ is regarded as a homology class
in the Borel-Moore homology $H^{BM}_{2*}(Y;\Z)$.
\end{cor}
By letting $\alpha=1_X$ in the formulae (\ref{eq4}) and resp.
(\ref{eq5}) above, we obtained the stratified multiplicative
property for the topological Euler-Poincar\'e characteristic and
resp. for the Chern-MacPherson class:
\begin{prop}\label{three} Let $f:X \to Y$ be a proper complex
algebraic (resp. analytic) map, with $Y$ irreducible (and compact in
the analytic context) and endowed with a complex algebraic (or
analytic) Whitney stratification $\VV$. Assume $f_*(1_X) \in
F_{\VV}(Y)$. Then:
\begin{equation}\label{eq6}
\chi(X)=\chi(F) \cdot \chi(Y) + \sum_{V < S} \left(
\chi({F_V})-\chi(F) \right) \cdot \hat{\chi}({\bar V}).
\end{equation}
\begin{equation}\label{eq7}
f_*(c_*(X))=\chi(F) \cdot c_*(Y) + \sum_{V < S} \left(
\chi({F_V})-\chi(F) \right) \cdot \hat{c}_*({\bar V}).
\end{equation}

\end{prop}

\section{Intersection homology Euler characteristics}

Let $Y$ be a topological pseudomanifold (or a locally cone-like
stratified space, [\cite{Sc}, p.232]), with a stratification $\VV$
by finitely many oriented strata of \emph{even} dimension. By
definition, strata of $\VV$ satisfy the frontier condition, and
$\VV$ is locally topologically trivial along each stratum $V$, with
fibers the cone on a compact pseudomanifold $L_{V,Y}$, the ``link"
of $V$ in $Y$. Note that each stratum $V$, its closure $\bar V$, and
in general any locally closed union of strata gets an induced
stratification of the same type. Examples are given by a complex
algebraic (or analytic) Whitney stratification of a reduced complex
algebraic (or compact complex analytic) variety.

Let $Sh_{\VV}(Y)$  be the category of $\VV$-constructible sheaves of
rational vector spaces, i.e., sheaves $\FC$ with the property that
for all $V \in \VV$ the restriction $\FC|_V$ is a locally constant
sheaf of $\Q$-vector spaces, with finite dimensional stalks. Denote
by $D^b_{\VV} (Y)$ the corresponding derived category of bounded
complexes with $\VV$-constructible cohomology sheaves (compare
\cite{B,KS,Sc}). Then one has an equality of Grothendieck groups
(e.g. compare [\cite{KS}, p. 77], [\cite{Sc}, Lemma 3.3.1])
$$K_0(Sh_{\VV}(Y)) = K_0(D^b_{\VV} (Y))$$ obtained by identifying the class
of a complex with the alternating sum of the classes of its
cohomology sheaves. Moreover one has a canonical group epimorphism
$$\chi_Y: K_0(D^b_{\VV} (Y)) \to F_{\VV}(Y)$$ defined by taking
stalkwise the Euler characteristic. Note that $\chi_Y$ is not
injective in general (e.g., see [\cite{Di}, p.98]), except for when
all strata $V \in \VV$ are simply-connected, e.g. for $Y = \{pt\}$,
in which case we use the shorter notion $K_0(pt)$. So $K_0(pt)$ is
just the Grothendieck group of finite dimensional $\Q$-vector
spaces, and it is a commutative ring with respect to tensor product,
with unit $\Q_{pt}$. Moreover, there is an isomorphism $K_0(pt)
\cong \Z$ induced by the Euler characteristic homomorphism.
$K_0(D^b_{\VV} (Y))$ becomes a unitary $K_0(pt)$-module, with the
multiplication defined by the exterior product: $$K_0(D^b_{\VV} (Y))
\times K_0(pt) \to K_0(D^b_{\VV \times \{pt\}} (Y \times
\{pt\}))=K_0(D^b_{\VV} (Y)),$$ and the Euler characteristic
homomorphisms $\chi_Y$ and $\chi$ are compatible with this structure
(more generally $\chi_Y$ commutes with exterior products).

Important examples of $\VV$-constructible complexes are provided by
the intersection cohomology complexes $IC_{\bar V}$ of the closures
of the strata $V \in \VV$, extended by $0$ to all of $Y$ (cf.
\cite{BBD,B,GM2}). These are selfdual with respect to Verdier
duality (and become important in the context of perverse sheaves and
mixed Hodge modules, as in our forthcoming papers \cite{CMS} and
\cite{CLMS}). The normalization axiom for $IC_{\bar V}$ (in the
conventions of \cite{BBD}) yields that $IC_{\bar V}|_V = \Q_V
[\text{dim}(V)]$, with $\text{dim}(V):= \text{dim}_{\R}(V)/2$ (the
complex dimension in the complex algebraic/analytic context). Since
we work in Grothendieck groups, in order to avoid signs in our
calculations, we will use the normalization condition of \cite{B},
that is, we work with $IC'_{\bar V}:=IC_{\bar V}[-\text{dim}(V)]$,
whose hypercohomology is exactly the intersection cohomology of
$\bar V$.

Let us fix for each $W \in \VV$ a point $w \in W$ with inclusion
$i_w:\{w\} \hookrightarrow Y$. Then \begin{equation}\label{100}
i_w^*[IC'_{\bar W}]=[i_w^*IC'_{\bar W}]=[\Q_{pt}]\in
K_0(w)=K_0(pt),\end{equation} and $i_w^*[IC'_{\bar V}] \neq [0] \in
K_0(pt)$ only if $W \leq V$. If we let
\begin{equation}\label{101} ic_{\bar V}:=\chi_Y(IC'_{\bar V}) \in
F_{\VV}(Y)\end{equation} be the corresponding constructible
function, then \begin{equation}\label{102} supp(ic_{\bar V})= \bar V
\ \text{ and } \ ic_{\bar V}|_V=1_V.\end{equation} Note that
$ic_{\bar V}(w)$ does not depend on the choice of $w \in W$, and
this is also the case for $i_w^*[IC_{\bar V}] \in K_0(pt)$. In fact,
since for any $j \in \Z$,
$$\HC^j (i_w^*IC'_{\bar V}) \simeq IH^j(c^{\circ}L_{W,V})$$ with
$c^{\circ}L_{W,V}$ the open cone on the link $L_{W,V}$ of $W$ in
$\bar V$ for $W \leq V$ (cf. \cite{B}, p.30, Prop.4.2), we have that
$$i_w^*[IC'_{\bar V}]=[IH^*(c^{\circ}L_{W,V})] \in K_0(pt).$$ In
terms of constructible functions, this gives
\begin{equation}\label{103} ic_{\bar V}(w)=I\chi(c^{\circ}L_{W,V}):=\chi([IH^*(c^{\circ}L_{W,V})]).\end{equation} In
particular, $\{ic_{\bar V} | V \in \VV \}$ is another distinguished
basis of $F_{\VV}(Y)$ since, by (\ref{102}), the transition matrix
to the basis $\{1_V\}$ is upper triangular with respect to $\leq$,
with all diagonal entries equal to $1$. Moreover, by (\ref{100}) the
$K_0(pt)$-submodule $\langle [IC'_{\bar V}] \rangle$ of
$K_0(D^b_{\VV} (Y))$ generated by the elements $[IC'_{\bar V}]$ ($V
\in \VV$) is in fact freely generated by them, and the restriction
\begin{equation}\label{iso} \chi_Y: \langle [IC'_{\bar V}] \rangle \to F_{\VV}(Y)
\end{equation}
is an isomorphism.

The main technical result of this section is the following
\begin{thm}\label{four} Assume $Y$ has an open dense stratum $S \in
\VV$ so that $V \leq S$ for all $V$. For each $V \in \VV \setminus
\{S\}$ define inductively
\begin{equation}\label{eq8}
\widehat{IC}(\bar V):=[IC'_{\bar V}] - \sum_{W < V}
\widehat{IC}(\bar W) \cdot i_w^* [IC'_{\bar V}] \in
K_0(D^b_{\VV}(Y)),
\end{equation}
and similarly
\begin{equation}\label{eq9}
\widehat{ic}(\bar V):=ic_{\bar V} - \sum_{W < V} \widehat{ic}(\bar
W) \cdot I\chi(c^{\circ}L_{W,V}) \in F_{\VV}(Y)
\end{equation}
so that $\chi_Y(\widehat{IC}(\bar V))=\widehat{ic}(\bar V)$. As the
notation suggests, $\widehat{IC}(\bar V)$ and $\widehat{ic}(\bar V)$
depend only on the stratified space $\bar V$ with its induced
stratification.
\begin{enumerate}
\item Assume $[\FC] \in K_0(D^b_{\VV}(Y))$ is an element of the
$K_0(pt)$-submodule $\langle [IC'_{\bar V}] \rangle$. Then
\begin{equation}\label{eq10}
[\FC]= [IC'_Y] \cdot i_s^*[\FC]+\sum_{V < S}  \widehat{IC}(\bar V)
\cdot \left( i_v^*[\FC] -i_s^*[\FC] \cdot i_v^*[IC'_Y] \right) \in
K_0(D^b_{\VV}(Y)).
\end{equation}
\item For any $\VV$-constructible function $\alpha \in F_{\VV}(Y)$, one has the equality
\begin{equation}\label{eq11}
\alpha=\alpha(s) \cdot ic_Y + \sum_{V < S} \left( \alpha(v) -
\alpha(s) \cdot I\chi(c^{\circ}L_{V,Y}) \right) \cdot
\widehat{ic}(\bar V).
\end{equation}

\end{enumerate}
\end{thm}

\begin{proof}
Note that the equation (\ref{eq11}) of the second part of the
theorem is a direct consequence of formula (\ref{eq10}) from the
first part. Indeed, by (\ref{iso}) we can first represent any
$\alpha \in F_{\VV}(Y)$ as $\alpha=\chi_Y([\FC])$, for some
$[\FC]\in \langle [IC'_{\bar V}] \rangle$. Then, assuming
(\ref{eq10}) holds for this choice of $[\FC]$, we can apply $\chi_Y$
to this equation and obtain (\ref{eq11}).

In order to prove formula (\ref{eq10}), consider
\begin{equation}\label{coef}[\FC]=\sum_V [IC'_{\bar V}] \cdot L(V),\end{equation}
for some $L(V) \in K_0(pt)$. The aim is to identify these
coefficients $L(V)$. Since $S$ is an open stratum, by applying
$i_s^*$ to (\ref{coef}) we obtain:
$$i_s^*[\FC]=L(S) \in K_0(s)=K_0(pt).$$
Next fix a stratum $W \neq S$, and apply $i_w^*$ to (\ref{coef}).
Recall that $i_w^*[IC'_{\bar W}]=[\Q_{pt}]\in K_0(w)=K_0(pt),$ and
$i_w^*[IC'_{\bar V}] \neq [0] \in K_0(pt)$ only if $W \leq V$. We
obtain \begin{equation}\label{1} i_w^*[\FC]=L(W)+\sum_{W < V}
i_w^*[IC'_{\bar V}] \cdot L(V) \in K_0(w)=K_0(pt).\end{equation}
Since $S$ is dense, we have that $W < S$, so the stratum $S$ appears
in the summation on the right hand side of (\ref{1}). Therefore
\begin{equation}\label{2}
i_w^*[\FC]-i_w^*[IC'_{Y}] \cdot i_s^*[\FC]=L(W)+\sum_{W<V<S}
i_w^*[IC'_{\bar V}] \cdot L(V) \in K_0(w)=K_0(pt).
\end{equation}
This implies that we can inductively calculate $L(V)$ in terms of
$$L'(W):=i_w^*[\FC]-i_w^*[IC'_{Y}] \cdot i_s^*[\FC].$$
Indeed, (\ref{2}) can be rewritten as \begin{equation}\label{3}
L'(W)=\sum_{W\leq V<S} i_w^*[IC'_{\bar V}] \cdot L(V) \in K_0(pt),
\end{equation}
and the matrix $A=(a_{W,V})$, with $a_{W,V}:=i_w^*[IC'_{\bar V}] \in
K_0(pt)$ for $W,V \in \VV \setminus \{S\}$, is upper-triangular with
respect to $\leq$, with ones on the diagonal. So $A$ can be
inverted. The non-zero coefficients of $A^{-1}=(a'_{W,V})$ can
inductively be calculated by $a'_{V,V}=1$ and
\begin{equation}\label{4}
a'_{W,V}=-\sum_{W \leq T < V} a'_{W,T} \cdot a_{T,V}
\end{equation}
for $W<V$. Then (\ref{coef}) becomes
\begin{equation}\label{5}[\FC]=[IC'_Y] \cdot i_s^*[\FC] + \sum_{W < S}  [IC'_{\bar W}] \cdot L(W)=
[IC'_Y] \cdot i_s^*[\FC] + \sum_{W \leq V< S}  [IC'_{\bar W}] \cdot
a'_{W,V} \cdot L'(V).
\end{equation}
The result follows by the inductive identification (for $V<S$
fixed):
$$\sum_{W\leq V} [IC'_{\bar W}] \cdot
a'_{W,V}=[IC'_{\bar V}]-\sum_{W \leq T <V} [IC'_{\bar W}] \cdot
a'_{W,T} \cdot a_{T,V}=\widehat{IC}(\bar V).$$

\end{proof}

\begin{remark}\rm
In this paper, we only make use of the equation (\ref{eq11}), and
this could be proven directly by working in $F_{\VV}(Y)$, following
the same arguments as above. However, the formula of equation
(\ref{eq10}) is particularly important since in the complex
algebraic context it extends to the framework of Grothendieck groups
of algebraic mixed Hodge modules that will be used in our
forthcoming paper \cite{CMS}. Of course, the technical condition
used in proving formula (\ref{eq10}) is not generally satisfied, but
it holds under the assumption of \emph{trivial monodromy along all
strata $V \in \VV$} (e.g., if all strata $V$ are simply-connected).
For more details, see \cite{CMS}.
\end{remark}

\bigskip

For the remaining part of this section, we will specialize to the
complex algebraic (or compact complex analytic) context, that is,
$Y$ is a reduced complex algebraic variety (or a reduced compact
complex analytic space), with a complex algebraic (resp. analytic)
Whitney stratification $\VV$. In this setting, let $f:X \to Y$ be a
proper complex algebraic (or analytic) map. Assume $f_*(\alpha) \in
F_{\VV}(Y)$ for a given $\alpha \in F_c(X)$, e.g., we choose $\VV'$
and $\VV$ complex Whitney stratifications of $X$ and resp. $Y$ such
that $f$ is a stratified submersion, and $\alpha \in F_{\VV'}(X)$.
Then
$$f_*(\alpha)=\sum_{V \in \VV} f_*(\alpha)(V) \cdot 1_V,$$ with
$$f_*(\alpha)(V)=\chi(\alpha|_{F_V}),$$ for $F_V$ the fiber of $f$
over a point in $V$. Assume, moreover, that $Y$ is irreducible, so
there is a dense stratum $S \in \VV$ with $F:=F_S$ a general fiber
of $f$. In terms of $f_*(\alpha)$, the equation (\ref{eq11}) of
Theorem \ref{four} becomes
\begin{equation}\label{eq12}
f_*(\alpha)=\chi(\alpha|_F) \cdot ic_Y + \sum_{V < S} \left(
\chi(\alpha|_{F_V})-\chi(\alpha|_F) \cdot I\chi(c^{\circ}L_{V,Y})
\right) \cdot \widehat{ic}({\bar V}) \in F_{\VV}(Y).
\end{equation}

By applying the group homomorphism $\chi$ and resp. $c_*$ to the
equation (\ref{eq12}), we obtain the following (recall $\chi \circ
f_*=\chi$ and $c_* \circ f_*=f_* \circ c_*$ for $f$ proper):
\begin{cor}\label{five}
\begin{equation}\label{eq13}
\chi(\alpha)=\chi(\alpha|_F) \cdot I\chi(Y) + \sum_{V < S} \left(
\chi(\alpha|_{F_V})-\chi(\alpha|_F) \cdot I\chi(c^{\circ}L_{V,Y})
\right) \cdot \widehat{I\chi}({\bar V}).
\end{equation}
\begin{equation}\label{eq14}
f_*(c_*(\alpha))=\chi(\alpha|_F) \cdot Ic_*(Y) + \sum_{V < S} \left(
\chi(\alpha|_{F_V})-\chi(\alpha|_F) \cdot I\chi(c^{\circ}L_{V,Y})
\right) \cdot \widehat{Ic}_*({\bar V}).
\end{equation}
Here $I\chi(Y):=\chi(ic_Y)=\chi([IH^*(Y;\Q)])$ is the intersection
homology Euler characteristic of $Y$, and similarly for
$\widehat{I\chi}({\bar V})$. Also $Ic_*(Y):=c_*(ic_Y)$, and
similarly for $\widehat{Ic}_*({\bar V})$, which by functoriality is
regarded as a homology class in $H^{BM}_{2*}(Y;\Z)$, the even degree
Borel-Moore homology of $Y$.
\end{cor}

By letting $\alpha=1_X$ in the formulae (\ref{eq13}) and resp.
(\ref{eq14}) above, we obtain:
\begin{prop}\label{six} Let $f:X \to Y$ be a proper complex
algebraic (resp. analytic) map, with $Y$ irreducible (and compact in
the analytic context) and endowed with a complex algebraic (or
analytic) Whitney stratification $\VV$. Assume $f_*(1_X) \in
F_{\VV}(Y)$. Then:
\begin{equation}\label{eq15}
\chi(X)=\chi(F) \cdot I\chi(Y) + \sum_{V < S} \left(
\chi({F_V})-\chi(F) \cdot I\chi(c^{\circ}L_{V,Y}) \right) \cdot
\widehat{I\chi}({\bar V}).
\end{equation}
\begin{equation}\label{eq16}
f_*(c_*(X))=\chi(F) \cdot Ic_*(Y) + \sum_{V < S} \left(
\chi({F_V})-\chi(F) \cdot I\chi(c^{\circ}L_{V,Y}) \right) \cdot
\widehat{Ic}_*({\bar V}).
\end{equation}

\end{prop}

In the special case when $f$ is the identity map, the equation
(\ref{eq15}) yields a formula expressing the Euler characteristics
of usual and intersection homology of an algebraic (or analytic)
variety in terms of each other and corresponding invariants of the
subvarieties formed by the closures of its singular strata.
Similarly, (\ref{eq16}) yields in this case a comparison between the
corresponding Chern homology classes of MacPherson:
\begin{cor} Let $Y$ be an irreducible complex algebraic (or compact analytic) variety endowed
with a complex algebraic (or analytic) Whitney stratification $\VV$.
Then in the above notations we have:
\begin{equation}\label{c1}
\chi(Y)=I\chi(Y) + \sum_{V < S} \left( 1- I\chi(c^{\circ}L_{V,Y})
\right) \cdot \widehat{I\chi}({\bar V}).
\end{equation}
\begin{equation}\label{c2}
c_*(Y)=Ic_*(Y) + \sum_{V < S} \left( 1- I\chi(c^{\circ}L_{V,Y})
\right) \cdot \widehat{Ic}_*({\bar V}).
\end{equation}
\end{cor}

The stratified multiplicative property for the intersection homology
Euler characteristic and for the corresponding homology
characteristic classes is obtained from (\ref{eq13}) and resp.
(\ref{eq14}) above in the case when $\alpha=ic_X$. Indeed, we have:
\begin{prop}\label{seven} Let $f:X \to Y$ be a proper complex
algebraic (resp. analytic) map, with $X$ pure dimensional and $Y$
irreducible (and compact in the analytic context). Assume $Y$ is
endowed with a complex algebraic (or analytic) Whitney
stratification $\VV$ so that $f_*(ic_X) \in F_{\VV}(Y)$. Then:
\begin{equation}\label{eq17}
I\chi(X)=I\chi(F) \cdot I\chi(Y) + \sum_{V < S} \left(
I\chi(f^{-1}(c^{\circ}L_{V,Y}))-I\chi(F) \cdot
I\chi(c^{\circ}L_{V,Y}) \right) \cdot \widehat{I\chi}({\bar V}).
\end{equation}
\begin{equation}\label{eq18}
f_*(Ic_*(X))=I\chi(F) \cdot Ic_*(Y) + \sum_{V < S} \left(
I\chi(f^{-1}(c^{\circ}L_{V,Y}))-I\chi(F) \cdot
I\chi(c^{\circ}L_{V,Y}) \right) \cdot \widehat{Ic}_*({\bar V}).
\end{equation}

\end{prop}

\begin{proof} Based on the above considerations, it suffices to show
that
\begin{equation}\label{eq19}\chi(ic_X|_F)=I\chi(F)\end{equation} and
\begin{equation}\label{eq20}\chi(ic_X|_{F_V})=I\chi(f^{-1}(c^{\circ}L_{V,Y})).\end{equation}

Since the general fiber $F$ of $f$ is locally normally nonsingular
embedded in $X$, we have a quasi-isomorphism (\cite{GM2}, \S 5.4.1):
$$IC'_X |_{F} \simeq IC'_{F},$$
hence an equality $ic_X|_F=ic_F$, thus proving (\ref{eq19}).

Similarly, since $\chi(ic_X|_{F_V})=f_*(ic_X)(v)$, for some $v \in
V$, in order to prove (\ref{eq20}), it suffices to show that
\begin{equation}\label{eq21} \HC^j(Rf_*IC'_X)_{v} \cong
IH^{j}(f^{-1}(c^\circ L_{V,Y});\Q).\end{equation} Let $N$ be a
normal slice to $V$ at $v$ in local analytic coordinates $(Y,v)
\hookrightarrow (\C^n,v)$, that is, a germ of a complex manifold
$(N,v) \hookrightarrow (\C^n,v)$, intersecting $V$ transversally
only at $v$, and with ${\rm dim} V + {\rm dim} N=n$. Recall that the
link $L_{V,Y}$ of the stratum $V$ in $Y$ is defined as
$$L_{V,Y}:=Y \cap N \cap
\partial B_r(v),$$
where $B_r(v)$ is an open ball of (very small) radius $r$ around
$v$. Moreover, $Y \cap N \cap B_r(v)$ is isomorphic (in a stratified
sense) to the open cone $c^{\circ}L_{V,Y}$ on the link (\cite{B}, p.
44). By factoring $i_v$ as the composition:
$$\{v\} \overset{\phi}{\hookrightarrow} Y \cap N
\overset{\psi}{\hookrightarrow} Y$$ we can now write:
{\allowdisplaybreaks
\begin{eqnarray*}
\HC^j(Rf_*IC'_X)_{v}  &\cong& \HB^j(Y,{i_v}_*i_v^*Rf_*IC'_X)\\
&\cong& \HB^j(v,\phi^*\psi^*Rf_*IC'_X)\\
&\cong& \HC^j(\psi^*Rf_*IC'_X)_{v}\\
&\cong& \HB^j(c^{\circ}L_{V,Y},Rf_*IC'_X)\\
&\cong& \HB^j(f^{-1}(c^{\circ}L_{V,Y}),IC'_X)\\
&\overset{(1)}{\cong}& \HB^j(f^{-1}(c^{\circ}L_{V,Y}),IC'_{f^{-1}(c^{\circ}L_{V,Y})})\\
&\cong& IH^{j}(f^{-1}(c^{\circ}L_{V,Y});\Q)
\end{eqnarray*}
} where in $(1)$ we used the fact that the inverse image of a normal
slice to a stratum of $Y$ in a stratification of $f$ is (locally)
normally non-singular embedded in $X$ (this fact is a consequence of
first isotopy lemma).

\end{proof}

\providecommand{\bysame}{\leavevmode\hbox
to3em{\hrulefill}\thinspace}


\begin{thebibliography}{10}


\bibitem{BBD}  Beilinson, A. A.,  Bernstein, J.,  Deligne, P.,  \emph{Faisceaux pervers},  Analysis and topology on singular spaces, I (Luminy, 1981),
5--171, Ast\'erisque, 100, Soc. Math. France, Paris, 1982.

\bibitem{B} Borel, A. et. al. , \emph{Intersection cohomology}, Progress in Mathematics,
vol. 50, Birkh\"auser Boston, Boston, MA, 1984.

\bibitem{CS1} Cappell, S. E., Shaneson, J. L., \emph{Genera of algebraic varieties and counting of lattice points},  Bull. Amer. Math. Soc. (N.S.)  30  (1994),  no. 1, 62--69.

\bibitem{CMS0} Cappell, S. E., Maxim, L. G., Shaneson, J. L.,
\emph{Euler characteristics of algebraic varieties},
math.AT/0606654, v2.


\bibitem{CMS} Cappell, S. E., Maxim, L., Shaneson, J. L.,
\emph{Hodge genera of algebraic varieties, I.}, to appear in Comm.
in Pure and Applied Math.

\bibitem{CLMS} Cappell, S. E., Libgober, A., Maxim, L., Shaneson, J. L.,
\emph{Hodge genera of algebraic varieties, II.}, math.AG/0702380.


\bibitem{Cat} de Cataldo, M. A., Migliorini, L., \emph{The Hodge theory of algebraic maps},
Ann. Scient. Ec. Norm. Sup., 4e serie, t.38, 2005, p. 693-750.



\bibitem{Di} Dimca, A., \emph{Sheaves in Topology}, Universitext, Springer-Verlag, Berlin, 2004.

\bibitem{F} Fulton, W., \emph{Introduction to toric varieties},
Annals of Mathematics Studies, 131. The William H. Roever Lectures
in Geometry. Princeton University Press, Princeton, NJ, 1993.

\bibitem{GM2} Goresky, M., MacPherson, M., \emph{Intersection Homology II}, Invent. Math., 71 (1983), 77-129.

\bibitem{KS} Kashiwara, M., Schapira, P., \emph{Sheaves on
Manifolds}, Grundlehren der Mathematischen Wissenschaften, 292.
Springer-Verlag, Berlin, Heidelberg, 1990.

\bibitem{M} MacPherson, R., \emph{Chern classes for singular algebraic varieties},
Ann. of Math. (2)  100  (1974), 423--432.

\bibitem{S} Shaneson, S., \emph{Characteristic classes, lattice points and Euler-MacLaurin formulae},
Proceedings of the International Congress of Mathematicians, Vol. 1,
2 (Z\"urich, 1994),  612--624, Birkh\"auser, Basel, 1995.

\bibitem{Sc} Sch\"urmann, J., \emph{Topology of singular spaces and constructible sheaves}, Monografie
Matematyczne, 63. Birkh\"auser Verlag, Basel, 2003.

\bibitem{SY} Sch\"urmann, J., Yokura, S., \emph{A survey of characteristic classes of singular spaces},
in ``Singularity Theory" (ed. by D. Ch\'eniot et al), Dedicated to
Jean Paul Brasselet on his $60$th birthday, Proceedings of the 2005
Marseille Singularity School and Conference, World Scientific, 2007,
865-952.

\bibitem{Sta} Stanley, R.P., \emph{Enumerative combinatorics, Volume 1},
Cambridge Studies in Advanced Mathematics 49, Cambridge University
Press, 2002.

\end{thebibliography}
\end{document}